\documentclass[12pt]{article} 
\usepackage{amsmath,amssymb,latexsym,theorem,bbm,bm}
\newfont{\msa}{msam10 scaled\magstep1}
\newfont{\ssmsa}{msam9}
\newfont{\hcmr}{cmr17 scaled\magstep5}

\newcommand{\SC}{\scriptstyle}

\newcommand{\CC}{\mathbb{C}}

\newcommand{\EE}{\mathsf{E}}
\newcommand{\csD}{{\SC\cal D}}

\newcommand{\PP}{\mathsf{P}}
\newcommand{\NN}{\mathbb{N}}
\newcommand{\RR}{\mathbb{R}}
\newcommand{\ZZ}{\mathbb{Z}}
\newcommand{\Bi}{\mathrm{Bi}}
\newcommand{\Be}{\mathrm{Be}}
\newcommand{\CP}{\mathrm{CP}}
\newcommand{\Po}{\mathrm{Po}}
\newcommand{\ee}{\mathrm{e}}
\newcommand{\dd}{\mathrm{d}}

\newcommand{\cL}{\mathcal{L}}

\newcommand{\tF}{\widetilde{F}}

\newcommand{\tR}{\widetilde{R}}

\newcommand{\distr}{\stackrel{\csD}{\longrightarrow}}

\newcommand{\sleq}{\mbox{\ssmsa\hspace*{0.1mm}\symbol{54}\hspace*{0.1mm}}}
\newcommand{\sgeq}{\mbox{\ssmsa\hspace*{0.1mm}\symbol{62}\hspace*{0.1mm}}}
\renewcommand{\leq}{\mbox{\msa\hspace*{0.9mm}\symbol{54}\hspace*{0.9mm}}}
\renewcommand{\geq}{\mbox{\msa\hspace*{0.9mm}\symbol{62}\hspace*{0.9mm}}}

\newcommand{\vare}{\varepsilon}
\newcommand{\varr}{\varrho}

\newcommand{\proofend}{\hfill$\square$}

\newcommand{\conv}{\operatornamewithlimits{\ast}}

\renewcommand{\Re}{\operatorname{Re}}

\numberwithin{equation}{section}

\theorembodyfont{\em}
\newtheorem{Lem}{Lemma}
\newtheorem{Thm}{Theorem}

\theorembodyfont{\rm}
\newtheorem{Rem}{Remark}
\newtheorem{Ex}{Example}

\begin{document}

\bibliographystyle{plain}

\title{
Poisson limit of an inhomogeneous nearly critical INAR(1) model
\thanks{
The first author acknowledges the support of the Computer and
Automation Research Institute of the Hungarian Academy of Sciences.
The second and third authors have been supported by the Hungarian Scientific
 Research Fund under Grant No.\ OTKA-T048544/2005.
}}

\author{
\textsc{L\'aszl\'o Gy\"orfi, M\'arton Isp\'any, Gyula Pap} and
\textsc{Katalin Varga} \\[6mm]
        \textit{Cs\"org\H{o} S\'andor professzor hatvanadik
                sz\"ulet\'esnapj\'ara,} \\  
        \textit{tisztelettel, bar\'ats\'aggal}
}

\maketitle

\begin{abstract}
An inhomogeneous first--order integer--valued autoregressive (INAR(1)) process
is investigated, where the autoregressive type coefficient slowly converges to 
one. It is shown that the process converges weakly to a Poisson or a compound
Poisson distribution.

\emph{Keywords:} Integer--valued time series; INAR(1) model; nearly unstable
          model; Poisson approximation; Galton--Watson process

\emph{AMS 2000 Subject Classification:} Primary 60J80, Secondary 60J27; 60J85 
\end{abstract}

\section{Introduction}
\label{intr}

A zero start inhomogeneous first order integer-valued autoregressive (INAR(1))
 time series \ $(X_n)_{n\in\ZZ_+}$ \ is defined as 
 \begin{equation}\label{INAR1}
   \begin{cases}
    X_n = \sum\limits_{j=1}^{X_{n-1}} \xi_{n,j} + \vare_n,, & n \in \NN,\\
    X_0 = 0,
   \end{cases}
 \end{equation}
 where \ $\{ \xi_{n,j}, \, \vare_n : n, j \in \NN \}$ \ are independent
 non-negative integer-valued random variables such that
 \ $\{ \xi_{n,j} : j \in \NN \}$ \ are identically distributed and
 \ $\PP( \xi_{n,1} \in \{0,1\} ) = 1$ \ for each \ $n \in \NN$. 
\ In fact, \ $(X_n)_{n\in\ZZ_+}$ \ is a special Galton-Watson branching
 process with immigration such that the offspring distributions are Bernoulli
 distributions.
We can interpret \ $X_n$ \ as the size of the \ $n^\mathrm{th}$ \ generation of
 a population, \ $\xi_{n,j}$ \ is the number of offspring produced by the
 \ $j^\mathrm{th}$ \ individual belonging to the \ $(n-1)^\mathrm{th}$
 \ generation, and \ $\vare_n$ \ is the number of immigrants in the
 \ $n^\mathrm{th}$ \ generation.

The process \eqref{INAR1} is called INAR(1) since it may also be written in the
 form
 \[
   \begin{cases}
    X_n = \varr_n \circ X_{n-1} + \vare_n, & n\in\NN,\\
    X_0 = 0,
   \end{cases}
 \]
 where 
 \[
   \varr_n := \EE\xi_{n,1}
 \]
 denotes the mean of the Bernoulli offspring distribution in the
 \ $n^\mathrm{th}$ \ generation, and we use the
 \emph{Steutel and van Harn operator\/} \ $\varr \, \circ$ \ which is defined
 for \ $\varr \in [0,1]$ \ and for a non-negative integer-valued random
 variable \ $X$ \ by
 \[
   \varr \circ X
   :=\begin{cases}
      \sum\limits_{j=1}^X \xi_j, & X>0,\\[2mm]
      0, & X=0,
     \end{cases}
 \]
 where the counting sequence \ $(\xi_j)_{j\in\NN}$ \ consists of independent
 and identically distributed Bernoulli random variables with mean \ $\varr$,
 \ independent of \ $X$ \ (see Steutel and van Harn \cite{SH}), and the
 counting sequences involved in \ $\varr_n \circ X_{n-1}$, \ $n\in\NN$, \ are
 mutually independent and independent of \ $(\vare_n)_{n\in\NN}$.

Let us denote the factorial moments of the immigration distributions by
 \[
   m_{n,k} := \EE \vare_n ( \vare_n - 1 ) \cdots ( \vare_n - k + 1 ),
   \qquad n,k\in\NN.
 \]
If \ $m_{n,1} < \infty$ \ for all \ $n \in \NN$ \ then we have the recursion
 \[
   \EE X_n = \varr_n \EE X_{n-1} + m_{n,1}, \qquad n\in\NN,
 \]
 since
 \[
   \EE( X_n \mid X_{n-1} )
   = \EE\left( \sum_{j=1}^{X_{n-1}} \xi_{n,j}
               + \vare_n \,\bigg|\, X_{n-1} \right)
   = \sum_{j=1}^{X_{n-1}} \EE\xi_{n,j} + \EE \vare_n
   = X_{n-1} \varr_n + m_{n,1}.
 \]
Consequently, the sequence \ $(\varr_n)_{n\in\NN}$ \ of the offspring means
 plays a crucial role in the asymptotic behavior of the sequence
 \ $(X_n)_{n\in\ZZ_+}$ \ as \ $n\to\infty$.
\ The INAR(1) process \ $(X_n)_{n\in\ZZ_+}$ \ is called
 \emph{nearly critical\/} if \ $\varr_n\to 1$ \ as \ $n\to\infty$.
\ We will investigate the asymptotic behavior of nearly critical INAR(1)
 processes.

Non--negative integer--valued time series, known as counting processes, arise 
in several fields of medicine (see, e.g., Cardinal et al. \cite{CRL} and Franke 
and Seligmann \cite{FS}). To model counting processes Al--Osh and Alzaid 
\cite{AA1} proposed the INAR(1) model. Isp\'any et al. \cite{IPZ1} investigated
the asymptotic inference for nearly unstable INAR(1) models. Later on Al--Osh 
and Alzaid \cite{AA2} and Du and Li \cite{DL} generalized this model by 
introducing the INAR($p$) model. 

The INAR models are special branching processes where the offspring distributions
are Bernoulli distributions. The theory of branching processes has been developed
for a long time, see Athreya and Ney \cite{AN}, and it can be applied in various 
fields. Branching processes are well-known models of binary search trees, see
Devroye \cite{D}. A recent application of them is the domain of peer-to-peer file 
sharing networks. Traffic measurements show that the workload generated by P2P 
applications is the dominant part of most of the Internet segments. The file 
population dynamics can be described by these mathematical models which also make 
possible the design and control of peer-to-peer systems, see Adar and Huberman 
\cite{AH}, Zhao et al. \cite{ZSR}. Space-time processes are standard models in 
seismology, see Lise and Stella \cite{LS}. One of these is the Epidemic-type 
Aftershock Sequence (ETAS ) model and they serve for surveillance of infections 
diseases as well, see  Farrington et al. \cite{FKG}. The theory of branching 
processes can also be applied to data on different aspects of biodiversity or 
macroevolution by the help of using phylogenetic trees, see, e.g., Aldous and 
Popovic \cite{AP} and Haccou and Iwasa \cite{HI}. An inhomogeneous branching
mechanism has been considered in Isp\'any et al. \cite{IPZ2}. Drost et al.
\cite{DAW} proved that the limit experiment of a homogeneneous INAR(1) model
has a Poisson distribution.

The present paper seems to be the first attempt to deal with the so--called
nearly unstable inhomogeneous INAR(1) model. The paper is organized as follows.
In Section 2 two basic lemmas are proved for inhomogeneous INAR(1) process.
In Section 3 the case of Bernoulli immigrations, in Section 4 the case of 
non-Bernoulli immigrations with Poisson limit distribution are considered.
Section 5 is devoted to the general case when the limit distribution is a
compound Poisson distribution. The results are extended for triangular system 
of mixtures of binomial distributions. In the Appendix at the end of paper some
technical lemmas are gathered.

\section{Preliminaries}

Let \ $\Be(p)$ \ denote a Bernoulli distribution with mean \ $p \in [0,1]$.
\ The distribution of a random variable \ $\xi$ \ will be denoted by
 \ $\cL(\xi)$.
\ Consider the unit disk \ $D := \{z \in \CC : |z| \leq 1 \}$ \ of the complex
 plane \ $\CC$.
\ The (probability) generating function of a non-negative integer-valued random
 variable \ $\xi$ \ is given by \ $z \mapsto \EE (z^\xi)$ \ for \ $z \in D$,
 \ and we have \ $\EE (z^\xi) \in D$ \ for all \ $z \in D$.
\ Introduce the generating functions
 \[
   F_n(z) := \EE (z^{X_n}), \qquad
   G_n(z) := \EE (z^{\xi_{n,1}}), \qquad 
   H_n(z) := \EE (z^{\vare_n}),\qquad z\in D.
 \]

\begin{Lem}\label{Fn} \ 
For an arbitrary inhomogeneous INAR(1) process \ $(X_n)_{n\in\ZZ_+}$ \ we have
 \[
   F_n(z)
   = \prod_{k=1}^n H_k\big( 1 + \varrho_{[k,n]} (z-1) \big),
   \qquad n\in\NN,
 \]
 for all \ $z \in D$, \ where
 \[
   \varrho_{[k,n]} := \begin{cases}
                     \prod\limits_{\ell=k+1}^n \varr_\ell
                      & \text{for \ $1\leq k\leq n-1$,} \\
                     1 & \text{for \ $k=n$.}
           \end{cases}
 \]
\end{Lem}

\noindent
\textbf{Proof.} \
\ The basic recursion for the generating functions \ $F_n$, \ $n \in \NN$, \ is
 \begin{equation}\label{genrec}
  \begin{split}
   F_n(z) = & \EE \Big( z^{\sum_{j=1}^{X_{n-1}} \xi_{n,j} + \vare_n} \Big)
          = \EE \Big( \EE\Big( z^{\sum_{j=1}^{X_{n-1}}\xi_{n,j} + \vare_n}
                      \,\Big|\, X_{n-1}\Big) \Big)\\
          = & \EE\big(G_n(z)^{X_{n-1}}\big) H_n(z)
          = F_{n-1}(G_n(z)) H_n(z),
  \end{split}
 \end{equation}
 valid for all \ $z\in\CC$ \ with \ $z \in D$ \ and \ $G_n(z) \in D$, \ see
 Athreya and Ney \cite[p. 263]{AN}.
Clearly \ $z \in D$ \ implies \ $G_n(z) \in D$, \ hence \eqref{genrec} is valid
 for all \ $z \in D$.
\ Since \ $\cL(\xi_{n,1})=\Be(\varr_n)$, \ we have
 \[
   G_n(z) =  1-\varr_n + \varr_n z = 1 + \varr_n (z-1)      
 \]
 for all \ $z\in\CC$.
\ We prove the statement of the lemma by induction.
For \ $n=1$, \ we have \ $F_1(z)=H_1(z)=H_1(1+(z-1))$. 
\ By the recursion \eqref{genrec}, we obtain for \ $n\geq 2$
 \begin{align*}
  F_n(z) &= F_{n-1}\big( 1 + \varr_n (z-1) \big) H_n(z) 
          = H_n(z)
            \prod_{k=1}^{n-1}
             H_k\big( 1 + \varr_n \varrho_{[k,n-1]} (z-1) \big) \\
         &= \prod_{k=1}^n
             H_k\big( 1 + \varrho_{[k,n]} (z-1) \big),
 \end{align*}
 and the proof is complete.
\proofend 

In fact, \ $X_n$ \ can be considered as a sum of independent Galton-Watson
 processes without immigration. 
Namely,
 \begin{equation}\label{XYrep}
  X_n = \sum_{k=1}^n Y_{n,k}, \qquad n\in\NN,
 \end{equation}
 where
 \begin{equation}\label{Yrecursion}
  Y_{n,k} := \begin{cases}
             0 & \text{for \ $k=0$,} \\
             \sum\limits_{j=Y_{n-1,\,1}+\cdots+Y_{n-1,\,k-1}+1}^{Y_{n-1,\,1}+\cdots+Y_{n-1,\,k}}
              \xi_{n,j} & \text{for \ $1\leq k\leq n-1$,} \\
             \vare_n & \text{for \ $k=n$.}
            \end{cases}
 \end{equation}
The distribution of \ $Y_{n,k}$ \ is a mixture of binomial distributions with a
 common probability parameter \ $\varr_n$, \ since the number \ $Y_{n-1,k}$ \ of
 Bernoulli random variables in the sum \eqref{Yrecursion} is a random variable
 as well. 
For a probability measure \ $\mu$ \ on \ $\ZZ_+$ \ and for a number \ $p\in[0,1]$,
 \ the mixture \ $\Bi(\mu,p)$ \ of binomial distributions with parameters
 \ $\mu$ \ and \ $p$ \ is a probability measure \ on \ $\ZZ_+$ \ defined by
 \[
   \Bi(\mu,p)\{j\} := \sum_{\ell=j}^\infty 
                 \binom{\ell}{j} p^j (1-p)^{\ell-j} \mu \{\ell\} \qquad
   \text{for \ $j\in\ZZ_+$.}
 \]
It is a particular example for mixture of distributions, see Johnson and Kotz
 \cite[Section I.7.3]{JK}, because the common method for mixture of binomial
 distributions is to use different values of probability parameter, see Johnson
 and Kotz \cite[Section III.11]{JK}. 
Note that \ $\Bi(\mu,1) = \mu$.

\begin{Lem}\label{binomrep}
For all \ $n\in\NN$, \ $1\leq k\leq n$, \ the distribution of \ $Y_{n,k}$ \ is a
 mixture of binomial distributions with parameters \ $\vare_k$ \ and
 \ $\varrho_{[k,n]}$. 
\ Thus
 \[
   \cL(X_n) = \conv_{k=1}^n \Bi\big( \cL(\vare_k), \varrho_{[k,n]} \big), 
 \]
 where \ $\ast$ \ denotes convolution of probability measures.
\end{Lem}

\noindent
\textbf{Proof.} \
First we check that \ $\Bi\big( \Bi(\mu,p), q \big) = \Bi(\mu, pq)$ \ for an
 arbitrary probability measure \ $\mu$ \ on \ $\ZZ_+$ \ and for all
 \ $p,q\in[0,1]$.
\ Indeed, for all \ $j\in\ZZ_+$,
 \begin{align*}
  \Bi\big( \Bi(\mu, p), q \big) \{j\}
   &= \sum_{\ell=j}^\infty
       \binom{\ell}{j} q^j (1-q)^{\ell-j} \, \Bi(\mu, p) \{\ell\} \\ 
   &= \sum_{\ell=j}^\infty
       \binom{\ell}{j} q^j (1-q)^{\ell-j}
       \sum_{k=\ell}^\infty \binom{k}{\ell} p^\ell (1-p)^{k-\ell} \, \mu \{k\} \\
  &= \sum_{k=j}^\infty
      (pq)^j \, \mu \{k\}
      \sum_{\ell=j}^k
       \binom{\ell}{j} \binom{k}{\ell} (p(1-q))^{\ell-j} (1-p)^{k-\ell}\\
  &= \sum_{k=j}^\infty \binom{k}{j} (pq)^j \, \mu \{k\}
     \sum_{\ell=j}^k \binom{k-j}{k-\ell} (p(1-q))^{\ell-j} (1-p)^{k-\ell}\\
  &= \sum_{k=j}^\infty \binom{k}{j} (pq)^j (1-pq)^{k-j} \, \mu \{k\}\\
  &= \Bi(\mu, pq) \{j\}.
 \end{align*}
Since \ $Y_{1,1} = \vare_1$, \ thus \ $\cL(Y_{1,1}) = \Bi(\cL(\vare_1), 1)$, \ and
 \ $\cL(Y_{n,k}) = \Bi(\cL(Y_{n,k-1}), \varr_n)$ \ for all \ $n \geq 2$ \ and all
 \ $k = 1, \dots, n$, \ we obtain the statement of the lemma by induction using
 the previous argument.
\proofend 

Remark that Lemma \ref{binomrep} implies the formula given for the generating
 function of \ $X_n$ \ in Lemma \ref{Fn}, since the generating function of a
 distribution \ $\Bi(\mu, p)$ \ is \ $z \mapsto H( 1 + (z-1) p )$, \ where
 \ $H$ \ denotes the generating function of \ $\mu$. 

\section{Poisson limit distribution: the case of Bernoulli immigrations}

First consider the simplest case, when \ $\cL(\vare_n) = \Be(m_{n,1})$,
 \ $n \in \NN$.

\begin{Thm}\label{Bernoulli} \
Let \ $(X_n)_{n\in\ZZ_+}$ \ be an INAR(1) process such that
 \ $\PP( \vare_n \in \{0,1\} ) = 1$ \ for all \ $n \in \NN$.
\ Assume that
 \begin{enumerate}
  \item[\textup{(i)}] \ 
   $\varr_n < 1$ \ for all \ $n \in \NN$,
    \ $\lim\limits_{n\to \infty} \varr_n = 1$,
    \ $\sum\limits_{n=1}^{\infty} (1-\varr_n) = \infty$,
  \item[\textup{(ii)}] \
   $\lim\limits_{n\to\infty} \frac{m_{n,1}}{1-\varr_n} = \lambda \in [0, \infty)$.
 \end{enumerate}
Then 
 \begin{equation}\label{Polim}
   X_n \distr \Po(\lambda) \quad \text{as} \quad n\to\infty.
 \end{equation}
(Here and in the sequel \ $\Po(0)$ \ is understood as a Dirac measure
 concentrated at the point 1.)
\end{Thm}

Remark that the condition \ $\sum\limits_{n=1}^{\infty} (1-\varr_n) = \infty$
 \ may be replaced by \ $\prod\limits_{n=1}^{\infty} \varr_n = 0$.
\ Moreover, if \ $\lambda > 0$ \ then the condition
 \ $\lim\limits_{n\to \infty} \varr_n = 1$ \ may be replaced by
 \ $\lim\limits_{n\to \infty} m_{n,1} = 0$.

\noindent
\textbf{Proof.} \
In order to prove the statement, we will show that
 \[
   \lim_{n\to \infty} F_n(z) = \ee^{\lambda (z-1)}
 \]
 for all \ $z \in D$.
\ Since \ $\cL(\vare_j) = \Be(m_{j,1})$, \ we have
 \[
   H_k(z) = 1 + m_{k,1} (z-1),\qquad z \in D, \quad k \in \NN.      
 \]
Applying Lemma \ref{Fn}, we can write
 \begin{equation}\label{F_B}
  F_n(z) = \prod_{k=1}^n \left[ 1 + m_{k,1} \varrho_{[k,n]} (z-1) \right],\qquad
  z \in D, \quad n \in \NN.
 \end{equation}
Consider the functions \ $\tF_n:\CC\to\CC$, \ $n\in\NN$, \ defined by
 \begin{equation}\label{tF_B}
   \tF_n(z):= \prod_{k=1}^n \ee^{ m_{k,1} \varrho_{[k,n]} (z-1) }.
 \end{equation}
In fact, \eqref{tF_B} is the generating function of a Poisson distribution.
The terms in the products in \eqref{F_B} and \eqref{tF_B} are generating
functions of probability distributions, hence Lemma \ref{compprod} is 
applicable, and we obtain
 \[
   | \tF_n(z) - F_n(z) |
   \leq \sum_{k=1}^n
         \big| \ee^{ m_{k,1} \varrho_{[k,n]} (z-1) }
               - 1 - m_{k,1} \varrho_{[k,n]} (z-1) \big|
 \]
 for \ $z \in D$, \ $n \in \NN$.
\ An application of the inequality \ $|\ee^u -1-u| \leq |u|^2$ \ valid for all
 \ $u \in \CC$ \ with \ $|u| \leq 1/2$ \ implies
 \begin{equation}\label{Taylor_m1}
   \Big| \ee^{ m_{k,1} \varrho_{[k,n]} (z-1) }
         - 1 - m_{k,1} \varrho_{[k,n]} (z-1) \Big|
   \leq m_{k,1}^2 \varrho_{[k,n]}^2 |z-1|^2         
 \end{equation}
 for \ $z\in\CC$ \ with \ $m_{k,1} \varrho_{[k,n]} |z-1| \leq 1/2$.
\ By Lemma \ref{Toeplitz} and taking into account assumption
 \ $\lim\limits_{n\to\infty}\frac{m_{n,1}}{1-\varr_n}=\lambda\in[0,\infty)$, \ we
 have 
 \begin{equation}\label{max_m1}
  \max_{1\sleq k\sleq n} m_{k,1} \varrho_{[k,n]}
  = \max_{1\sleq k\sleq n} \frac{m_{k,1}}{1-\varr_k} a^{(1)}_{n,k}
  \to 0 \qquad
  \text{as \ $n\to\infty$.}     
 \end{equation}
Thus, the estimate \eqref{Taylor_m1} is valid for all \ $z \in D$, \ for
 sufficiently large \ $n$ \ and for all \ $k = 1, \dots, n$, \ and we obtain 
 \[
   | \tF_n(z) - F_n(z) |
   \leq  |z-1|^2 \sum_{k=1}^n m_{k,1}^2 \varrho_{[k,n]}^2 .    
 \]
By \ $\lim\limits_{n\to\infty} \frac{m^2_{n,1}}{1-\varr_n} = 0$ \ and by Lemma
 \ref{Toeplitz} we obtain
 \begin{equation}\label{quad}
   \sum_{k=1}^n m_{k,1}^2 \varrho_{[k,n]}^2
   = \sum_{k=1}^n \frac{m_{k,1}^2}{1-\varr_k} a^{(2)}_{n,k}
   \to0 \qquad \text{as \ $n\to\infty$.}
 \end{equation}
Consequently,
 \[
   \lim_{n\to\infty} | \tF_n(z) - F_n(z) | = 0 \qquad
   \text{for all \ $z \in D$.}
 \]
An application of Lemma \ref{Toeplitz} yields
 \begin{equation}\label{lim_m1}
   \sum_{k=1}^n m_{k,1} \varrho_{[k,n]}
   = \sum_{k=1}^n \frac{m_{k,1}}{1-\varr_k} a^{(1)}_{n,k}
   \to \lambda
   \qquad \text{as \ $n\to\infty$.}
 \end{equation}
Consequently,
 \[
   \lim_{n\to\infty} \tF_n(z) = \ee^{\lambda (z-1)} \qquad
   \text{for all \ $z \in D$,}
 \]
 and we obtain \ $F_n(z) \to \ee^{\lambda (z-1)}$ \ as \ $n \to \infty$ \ for all
 \ $z \in D$.
\proofend 

\noindent
\textbf{Second proof of Theorem \ref{Bernoulli} by Poisson approximation.} \
We may prove the theorem by Poisson approximation as well. 
The total variation distance between two probability measures \ $\mu$ \ and
 \ $\nu$ \ on \ $\ZZ_+$ \ equals
 \[
   d(\mu,\nu) = \frac{1}{2} \sum_{j=0}^\infty \big| \mu\{j\}-\nu\{j\} \big| .
 \]
A sequence \ $(\mu_n)_{n\in\NN}$ \ of probability measures on \ $\ZZ_+$
 \ converges weakly to a probability measure \ $\mu$ \ on \ $\ZZ_+$ \ if and
 only if \ $d(\mu_n,\mu)\to0$.
\ We prove \eqref{Polim} by showing that
 \begin{equation}\label{XPdist}
  d\big(\cL(X_n), \Po(\lambda)\big) \to 0 \qquad \text{as \ $n \to \infty$.}
 \end{equation}
One can easily check that \ $\Bi\big( \Be(p), q \big) = \Be(pq)$ \ for
 arbitrary \ $p,q\in[0,1]$, hence by Lemma \ref{binomrep} we obtain
 \ $\cL(X_n) = \conv\limits_{k=1}^n \Be\left( m_{k,1} \varrho_{[k,n]} \right)$.
\ By Lemma \ref{conv},
 \[
   d\Big(\cL(X_n),\,\conv_{k=1}^n \Po( m_{k,1} \varrho_{[k,n]} )\Big)
   \leq
   \sum_{k=1}^n
    d\big( \Be(m_{k,1} \varrho_{[k,n]}), \Po(m_{k,1} \varrho_{[k,n]}) \big).
 \]
We show that
 \begin{equation}\label{Be_Po}
  d\big( \Be(p), \Po(p) \big) \leq p^2
 \end{equation}
 for all \ $p\in [0,1]$. 
\ Indeed,
 \[
  d\big( \Be(p), \Po(p) \big)
  = \frac{1}{2} ( \ee^{-p} - 1 + p )
    + \frac{1}{2} ( p- p \, \ee^{-p} )
    + \frac{1}{2} ( 1 - \ee^{-p} - p \, \ee^{-p} )
  = p (1 - \ee^{-p} )
  \leq p^2.
 \]
Applying \eqref{Be_Po} and \eqref{quad}, we conclude
 \[
   d\Big(\cL(X_n),\,\conv_{k=1}^n \Po( m_{k,1} \varrho_{[k,n]} )\Big)
   \leq \sum_{k=1}^n m_{k,1}^2 \varrho_{[k,n]}^2
   \to 0.
 \]
Clearly,
 \[
   \conv_{k=1}^n \Po( m_{k,1} \varrho_{[k,n]} )
   = \Po\bigg( \sum_{k=1}^n m_{k,1} \varrho_{[k,n]} \bigg)
   \to \Po(\lambda)  
 \]
in law by \eqref{lim_m1}, and we obtain \ $X_n \distr \Po(\lambda)$.
\proofend 

\section{Poisson limit distribution: the case of non-Bernoulli immigrations}

\begin{Thm}\label{Poisson} \
Let \ $(X_n)_{n\in\ZZ_+}$ \ be an inhomogeneous INAR(1) process.
Assume that 
 \begin{enumerate}
  \item[\textup{(i)}] \ 
   $\varr_n < 1$ \ for all \ $n \in \NN$,
    \ $\lim\limits_{n\to \infty} \varr_n = 1$,
    \ $\sum\limits_{n=1}^{\infty} (1-\varr_n) = \infty$,
  \item[\textup{(ii)}] \
   $\lim\limits_{n\to\infty} \frac{m_{n,1}}{1-\varr_n} = \lambda \in [0, \infty)$,
   \ $\lim\limits_{n\to\infty} \frac{m_{n,2}}{1-\varr_n} = 0$.
 \end{enumerate}
Then 
 \[
   X_n \distr \Po(\lambda) \quad \text{as} \quad n\to\infty.
 \]
\end{Thm}

\begin{Rem}\label{domassump}
Since
 \[
   m_{n,1} = \sum_{j=1}^\infty \PP(\vare_n \geq j), \qquad
   m_{n,2} = 2 \sum_{j=1}^\infty j \PP(\vare_n > j),
 \]
 assumption (ii) implies
 \begin{equation}\label{taillim}
  \lim_{n\to\infty} \frac{\PP(\vare_n \geq 1)}{1-\varr_n} = \lambda, \qquad
  \lim_{n\to\infty} \frac{\PP(\vare_n \geq 2)} {1-\varr_n} = 0.
 \end{equation}
In general the converse is not true. 
However, if there exists a sequence \ $(b_j)_{j\in\NN}$ \ of non-negative real
 numbers such that \ $\sum_{j=1}^\infty j b_j < \infty$ \ and
 \ $\frac{\PP(\vare_n > j)}{1-\varr_n} \leq b_j$ \ for all \ $j, n \in \NN$,
 \ then \eqref{taillim} implies (ii) by the dominated convergence theorem.
\end{Rem}

\noindent
\textbf{Proof.} \
By Lemma \ref{Fn}, we can write
 \[
   F_n(z)
   = \prod_{k=1}^n H_k\big( 1 + \varrho_{[k,n]} (z-1) \big),
   \qquad z \in D, \quad n \in \NN,
 \]
 Consider the functions \ $\tF_n: \CC \to \CC$, \ $n \in \NN$, \ defined by
 \[
   \tF_n(z) = \prod_{k=1}^n \left[ 1 + m_{k,1} \varrho_{[k,n]} (z-1) \right].
 \]
By Lemma \ref{compprod}, we obtain
 \[
   | F_n(z) - \tF_n(z) |
   \leq \sum_{k=1}^n
         \Big| H_k\big( 1 + \varrho_{[k,n]} (z-1) \big)
               - 1 - m_{k,1} \, \varrho_{[k,n]} (z-1) \Big|
 \]
 for \ $z \in D$, \ $n \in \NN$.
\ Applying Lemma \ref{Taylor}, we have
 \[
   | H_k(u) - 1 - m_{k,1} (u-1) |
   \leq \frac{1}{2} m_{k,2} \, |u-1|^2 \qquad
   u \in D, \quad k \in \NN.
 \]
Thus
 \[
   \Big| H_k\big( 1 + \varrho_{[k,n]} (z-1) \big)
         - 1 - m_{k,1} \, \varrho_{[k,n]} (z-1) \Big|
   \leq \frac{1}{2} m_{k,2} \, \varrho_{[k,n]}^2 |z-1|^2
 \] 
 for all \ $z \in D$, \ since \ $z \in D$ \ implies
 \ $1 + \varrho_{[k,n]} (z-1) \in D$.
\ Consequently,
 \[
   | F_n(z) - \tF_n(z) |
   \leq \frac{1}{2} |z-1|^2 \sum_{k=1}^n m_{k,2} \, \varrho_{[k,n]}^2
   \to 0      
 \] 
as \ $n\to\infty$ \ for all \ $z \in D$ \ 
 by Lemma \ref{Toeplitz} taking into account assumption
 \ $\lim\limits_{n\to\infty} \frac{m_{n,2}}{1-\varr_n} = 0$.
\ Theorem \ref{Bernoulli} clearly implies \ $\tF_n(z) \to \ee^{\lambda(z-1)}$
 \ for all \ $z\in D$, \ hence we conclude \ $F_n(z)\to\ee^{\lambda (z-1)}$ \ as
 \ $n\to\infty$ \ for all \ $z\in D$.
\proofend 

\noindent
\textbf{Second proof of Theorem \ref{Poisson} by Poisson approximation.} \
Note that \ $m_{k,1} \varrho_{[k,n]} \leq 1$ \ for sufficiently large \ $n$ \ and
 for all \ $1\leq k\leq n$ \ by \eqref{max_m1}. 
By Lemmas \ref{binomrep} and \ref{conv}, we have, for sufficiently
 large \ $n\in\NN$,
 \[
   d\left(\cL(X_n), \, \conv_{k=1}^n \Be(m_{k,1} \varrho_{[k,n]}) \right)
   \leq \sum_{k=1}^n
          d\big( \Bi(\vare_k, \varrho_{[k,n]}),\,
                 \Be(m_{k,1} \varrho_{[k,n]}) \big).
 \]
We prove that
 \begin{equation}\label{tvarest}
  d\big(\Bi(\vare,p), \Be(p\EE\vare)\big)
  \leq \frac{3}{2} p^2 \EE\vare(\vare-1),
 \end{equation}
 where \ $p\in [0,1]$ \ and \ $\vare$ \ is a non-negative integer-valued 
 random variable such that \ $p\EE\vare\leq 1$. 
\ We have
 \[
   d\big(\Bi(\vare,p), \Be(p\EE\vare)\big)
   \leq \frac{1}{2}(A+B+C),
 \]
 where
 \begin{align*}
   A&:=\left|\sum_{\ell=0}^\infty (1-p)^\ell \, \PP(\vare=\ell)
              - (1 - p \EE\vare)\right|
    \leq \sum_{\ell=0}^\infty |(1-p)^\ell - 1 + \ell p| \, \PP(\vare=\ell), \\
   B&:=\left|\sum_{\ell=1}^\infty \ell p(1-p)^{\ell-1} \, \PP(\vare=\ell)
             -p \EE\vare \right|
    \leq p \sum_{\ell=1}^\infty \ell |(1-p)^{\ell-1} - 1| \, \PP(\vare=\ell),\\
   C&:=\sum_{j=2}^\infty \sum_{\ell=j}^\infty
        \binom{\ell}{j} p^j (1-p)^{\ell-j} \, \PP(\vare=\ell)
     = \sum_{\ell=2}^\infty \PP(\vare=\ell)
        \sum_{j=2}^\ell \binom{\ell}{j} p^j (1-p)^{\ell-j} \\
    &= \sum_{\ell=2}^\infty
        \PP(\vare=\ell) \bigl(1 - (1-p)^\ell - \ell p(1-p)^{\ell-1} \bigr)\\
    &\leq \sum_{\ell=2}^\infty
           \PP(\vare=\ell)
           \bigl( |1 - \ell p - (1-p)^\ell| + \ell p |1-(1-p)^{\ell-1}| \bigr).
 \end{align*}
By Taylor's formula for the function \ $p\mapsto (1-p)^k$ \ we get
 \[
   |1 - \ell p-(1-p)^\ell|
   \leq \frac{1}{2} \ell(\ell-1) p^2 \sup_{\theta\in[0,1]} (1-\theta p)^{\ell-2}
   \leq \frac{1}{2} \ell(\ell-1)p^2.
 \]
Finally, since
 \[
   |(1-p)^{\ell-1}-1|\leq p(\ell-1),
 \]
 we obtain \eqref{tvarest}. 
Thus, we have
 \[
   d\left(\cL(X_n), \, \conv_{k=1}^n \Be(m_{k,1} \varrho_{[k,n]}) \right)
   \leq \frac{3}{2} \sum_{k=1}^n m_{k,2} \, \varrho_{[k,n]}^2 ,
 \]
where the right hand side tends to 0 by the assumption
\ $\lim\limits_{n\to\infty} \frac{m_{n,2}}{1-\varr_n} = 0$. 
\ Obviously, Theorem \ref{Bernoulli} implies
 \[
   d\left( \, \conv_{k=1}^n \Be(m_{k,1} \varrho_{[k,n]}), \, \Po(\lambda)\right)
   \to 0 \qquad \text{as \ $n\to\infty$,}
 \]
 hence we obtain
 \[
   d\big(\cL(X_n), \Po(\lambda)\big) \to 0 \qquad \text{as \ $n\to \infty$,}
 \]
 which completes the proof.
\proofend 

In fact, a similar theorem holds for triangular system of mixtures of
binomial distributions.

\begin{Thm}\label{triangular}
Let \ $k_n\in\NN$ \ for all \ $n\in\NN$,\ and \ $\{\zeta_{n,k}:1\leq k\leq k_n,
\  n\in\NN\}$ \ be non-negative integer-valued random 
variables. Moreover, let \ $p_{n,k}\in [0,1]$, \ $1\leq k\leq k_n$,
\ $n\in\NN$. \ Assume that
 \begin{enumerate}
  \item[\textup{(i)}] \ 
  $\sum\limits_{k=1}^{k_n} p_{n,k}\EE\zeta_{n,k}\to \lambda$ \ for some
   \ $\lambda\geq 0$;
  \item[\textup{(ii)}] \
   $\sum\limits_{k=1}^{k_n} (p_{n,k}\EE\zeta_{n,k})^2\to 0$;
  \item[\textup{(iii)}] \  
   $\sum\limits_{k=1}^{k_n} p^2_{n,k}\EE\zeta_{n,k}(\zeta_{n,k}-1)\to 0$
 \end{enumerate}
as \ $n\to\infty$. \ Then
 \[
    \conv_{k=1}^{k_n} \Bi\big( \cL(\zeta_{n,k}), p_{n,k} \big)\to \Po(\lambda) 
 \]
in law as \ $n\to\infty$.
\end{Thm}

\section{Compound Poisson limit distribution}

Recall that if \ $\mu$ \ is a finite measure on \ $\ZZ_+$ \ then the
 \emph{compound Poisson distribution} \ $\CP(\mu)$ \ with intensity measure
 \ $\mu$ \ is the probability measure on \ $\ZZ_+$ \ with generating function
 \[
   z \mapsto \exp\left\{\sum_{j=1}^\infty \mu\{j\} (z^j-1)\right\}
   \qquad \text{for \ $z\in D$.}
 \]
In fact, \ $\CP(\mu)$ \ is an infinitely divisible distribution on \ $\ZZ_+$
 \ with L\'evy measure \ $\mu$ \ restricted onto \ $\NN$, \ and, for an
 arbitrary infinitely divisible distribution \ $\nu$ \ on \ $\ZZ_+$, \ there
 exists a finite measure \ $\mu$ \ on \ $\NN$ \ such that \ $\nu = \CP(\mu)$.
\ Moreover, \ $\CP(\mu)$ \ is the distribution of the random sum
 \[
   \sum_{j=1}^\Pi \zeta_j,
 \]
 where \ $\{\Pi, \, \zeta_j : j \in \NN\}$ \ are independent random variables,
 \ $\cL(\Pi) = \Po(\|\mu\|)$ \ and \ $\cL(\zeta_j) = \frac{\mu}{\|\mu\|}$ \ for
 \ $j \in \NN$, \ where \ $\|\mu\| := \sum_{j=1}^\infty \mu\{j\}$.
\ Further, \ $\CP(\mu)$ \ is the distribution of the weakly convergent infinite
 sum
 \[
   \sum_{j=1}^\infty j \eta_j,
 \]
 where \ $\{\eta_j : j\in\NN\}$ \ are independent random variables with
 \ $\cL(\eta_j) = \Po(\mu\{j\})$ \ for \ $j \in \NN$.
\ (See Barbour et al. \cite[Section 10.4]{BHJ}.)

First we consider the case when the intensity measure \ $\mu$ \ of the limiting
 compound Poisson distribution \ $\CP(\mu)$ \ has bounded support.

\begin{Thm}\label{CP_bounded} \
Let \ $(X_n)_{n\in\ZZ_+}$ \ be an inhomogeneous INAR(1) process.
Assume that
 \begin{enumerate}
  \item[\textup{(i)}] \ 
   $\varr_n < 1$ \ for all \ $n \in \NN$,
    \ $\lim\limits_{n\to \infty} \varr_n = 1$,
    \ $\sum\limits_{n=1}^{\infty} (1-\varr_n) = \infty$,
  \item[\textup{(ii)}] \
   $\lim\limits_{n\to\infty} \frac{m_{n,j}}{j(1-\varr_n)}
    = \lambda_j \in [0,\infty)$
   \ for \ $j = 1, \dots, J$ \ with \ $\lambda_J = 0$.   
 \end{enumerate}
Then
 \[
     X_n \distr \CP(\mu) 
   \quad \text{as} \quad n\to\infty,
 \]
 where \ $\mu$ \ is a finite measure on \ $\{1, \dots, J-1\}$ \ given by
 \begin{equation}\label{mudef}
  \mu\{j\} := \frac{1}{j!}
              \sum_{i=0}^{J-j-1} \frac{(-1)^i}{i!} \lambda_{j+i},
  \qquad j=1, \dots, J-1.
 \end{equation}
\end{Thm}

\begin{Rem}
One can easily check that \ $\lambda_j$, \ $j=1, \dots, J-1$, \ are the
first \ $J-1$ \ factorial moments of the measure \ $\mu$, \ i.e.,
\[
      \lambda_j = \sum_{i=j}^{J-1} i(i-1)\cdots(i-j+1) \mu\{i\}, 
      \qquad j=1, \dots, J-1.
\]
Moreover, since
\[
     m_{n,j} = j! \sum_{i=j}^\infty \binom{i-1}{j-1}\PP(\vare_n\geq i),
     \qquad j\in\NN,
\]
assumption (ii) implies
\begin{enumerate}
  \item[\textup{(ii)${}^\prime$}] \
   $\lim\limits_{n\to\infty} \frac{\PP(\vare_n\sgeq i)}{i(1-\varr_n)}
    = \mu \{i\} \in [0,\infty)$
   \ for \ $i = 1, \dots, J$ \ with \ $\mu\{J\} = 0$.   
 \end{enumerate}
On the other hand, (ii)${}^\prime$ and additional domination assumption,
see Remark \ref{domassump}, imply (ii). 
\end{Rem}

\noindent
\textbf{Proof.}
By Lemma \ref{Fn}, we can write
 \[
   F_n(z)
   = \prod_{k=1}^n H_k\big( 1 + \varrho_{[k,n]} (z-1) \big),
   \qquad z \in D, \quad n \in \NN.
 \]
Consider the functions
 \[
   \tF_n(z) := \prod_{k=1}^n \ee^{H_k( 1 + \varrho_{[k,n]} (z-1) ) - 1},
   \qquad z \in D, \quad n\in\NN.
 \]
By Lemma \ref{compprod}, we obtain
 \[
   | \tF_n(z) - F_n(z) |
   \leq \sum_{k=1}^n
         \Big| \ee^{H_k( 1 + \varrho_{[k,n]} (z-1) ) - 1}
               - H_k\big( 1 + \varrho_{[k,n]} (z-1) \big) \Big|
 \]
 for \ $z \in D$, \ $n \in \NN$.
\ An application of the inequality \ $|\ee^u -1-u| \leq |u|^2$ \ valid for all
 \ $u \in \CC$ \ with \ $|u| \leq 1/2$ \ implies
 \begin{equation}\label{Taylor_H}
  \Big| \ee^{H_k( 1 + \varrho_{[k,n]} (z-1) ) - 1}
        - H_k\big( 1 + \varrho_{[k,n]} (z-1) \big) \Big| 
  \leq \Big| H_k\big( 1 + \varrho_{[k,n]} (z-1) \big) - 1 \Big|^2       
 \end{equation}
 for \ $z\in\CC$ \ with
 \ $\Big| H_k\big( 1 + \varrho_{[k,n]} (z-1) \big) - 1 \Big| \leq 1/2$.
\ Applying Lemma \ref{Taylor}, we have
 \[
   | H_k(u) - 1 |
   \leq m_{k,1} \, |u-1|, \qquad
   u \in D, \quad k \in \NN.
 \]
Thus
 \[
   \Big| H_k\big( 1 + \varrho_{[k,n]} (z-1) \big) - 1 \Big|
   \leq m_{k,1} \, \varrho_{[k,n]} |z-1|
 \] 
 for all \ $z \in D$, \ since \ $z \in D$ \ implies
 \ $1 + \varrho_{[k,n]} (z-1) \in D$. \ By Lemma \ref{Toeplitz} and 
taking into account assumption \ $\lim\limits_{n\to\infty}\frac{m_{n,1}}
{1-\varr_n}=\lambda_1\in[0,\infty)$, \ we obtain \eqref{max_m1}.
Thus, the estimate \eqref{Taylor_H} is valid for all \ $z \in D$, \ for
 sufficiently large \ $n$ \ and for all \ $k = 1, \dots, n$, \ and we obtain 
 \[
   | \tF_n(z) - F_n(z) |
   \leq |z-1|^2 \sum_{k=1}^n m_{k,1}^2 \, \varrho_{[k,n]}^2
   \to 0 \qquad
   \text{as \ $n \to \infty$ \ for all \ $z \in D$}     
 \] 
 by \eqref{quad}.
Clearly
 \[
   \tF_n(z)
    = \exp\left\{ \sum_{k=1}^n
                   \Big[ H_k\big( 1 + \varrho_{[k,n]} (z-1) \big) - 1 \Big]
          \right\} .
 \]
Again by Lemma \ref{Taylor}, we have
 \[
   H_k\big( 1 + \varrho_{[k,n]} (z-1) \big) - 1
   = \sum_{j=1}^{J-1} \frac{m_{k,j}}{j!} \varr_{[k,n]}^j (z-1)^j + R_{n,k,J}(z),
 \]
 for all \ $z \in D$, \ for sufficiently large \ $n$ \ and for all
 \ $k = 1, \dots, n$, \ where
 \[
   | R_{n,k,J}(z) | \leq \frac{m_{k,J}}{J!} \varr_{[k,n]}^J |z-1|^J.
 \] 
An application of Lemma \ref{Toeplitz} yields
 \begin{equation}\label{lim_m}
   \sum_{k=1}^n m_{k,j} \varrho_{[k,n]}^j
   = \sum_{k=1}^n \frac{m_{k,j}}{1-\varr_k} a^{(j)}_{n,k}
   \to \lambda_j
   \qquad \text{as \ $n\to\infty$}
 \end{equation}
 for \ $j=1,\ldots,J$.
\ Moreover, by \eqref{lim_m},
 \[
   \sum_{k=1}^n | R_{n,k,J}(z) |
   \leq \frac{|z-1|^J}{J!} \sum_{k=1}^n m_{k,J} \varrho_{[k,n]}^J
   \to 0
 \]
as \ $n\to\infty$ \ for all \ $z\in D$ \ since \ $\lambda_J=0$. 
\ Consequently,
 \[
   \lim_{n\to\infty} F_n(z)
   = \exp\left\{ \sum_{j=1}^{J-1} \frac{\lambda_j}{j!} (z-1)^j \right\} 
   =: F(z) \qquad
   \text{for all \ $z \in D$,}
 \]
 where \ $F$ \ is the generating function of a probability distribution.
Clearly
 \begin{align*}
   \sum_{j=1}^{J-1} \frac{\lambda_j}{j!} (z-1)^j
   &= \sum_{j=1}^{J-1}
       \frac{\lambda_j}{j!}
       \sum_{i=0}^j \binom{j}{i} (-1)^{j-i} z^i
    = \sum_{i=0}^{J-1}
       \frac{z^i}{i!}
       \sum_{j=i}^{J-1}
        \frac{(-1)^{j-i}}{(j-i)!} \lambda_j \\
   &= \sum_{i=0}^{J-1}
       \frac{z^i}{i!}
       \sum_{j=0}^{J-i-1}
        \frac{(-1)^j}{j!} \lambda_{j+i}
    = \sum_{i=1}^{J-1}
       \mu\{i\} (z^i-1)
 \end{align*}
since \ $\sum_{i=1}^{J-1} \mu\{i\} = \sum_{j=1}^{J-1} 
\frac{(-1)^j}{j!} \lambda_j$, \
and we obtain \ $X_n \distr \CP(\mu)$ \ as \ $n \to \infty$.
\proofend 

\textbf{Second proof of Theorem \ref{CP_bounded} by Poisson approximation.}
 \
By Lemmas \ref{binomrep} and \ref{conv}, we have
 \[
   d\left(\cL(X_n), \,
          \conv_{k=1}^n \CP\big( \Bi(\vare_k, \varrho_{[k,n]}) \big) \right)
   \leq \sum_{k=1}^n
         d\Big( \Bi(\vare_k, \varrho_{[k,n]}), \,
                \CP\big( \Bi(\vare_k, \varrho_{[k,n]}) \big) \Big).
 \]
We prove that
 \begin{equation}\label{Bi_CP}
  d\Big(\Bi(\vare, p), \CP\big(\Bi(\vare, p)\big)\Big)
  \leq p^2 (\EE\vare)^2
 \end{equation}
 for all \ $p\in [0,1]$ \ and for all non-negative integer-valued 
 random variable \ $\vare$.
\ By Barbour et al. \cite[Corollary 10.L.1]{BHJ}, we have
 \[
   d\Big(\Bi(\vare, p), \CP\big(\Bi(\vare, p)\big)\Big)
   \leq \PP\big(\Bi(\vare, p) \geq 1\big)^2.
 \]
Now
 \begin{equation}\label{Bi}
  \begin{split}
   \PP(\Bi(\vare, p) \geq 1)
   & = 1 - \sum_{\ell=0}^\infty (1-p)^\ell \, \PP(\vare=\ell)
     = \sum_{\ell=0}^\infty \big[1 - (1-p)^\ell\big] \, \PP(\vare=\ell) \\
   & \leq \sum_{\ell=0}^\infty p \ell \, \PP(\vare=\ell)
     = p \, \EE\vare,
  \end{split}
 \end{equation}
 hence we obtain \eqref{Bi_CP}.
Applying \eqref{Bi_CP}, we conclude
 \[
   d\left(\cL(X_n), \,
          \conv_{k=1}^n \CP\big( \Bi(\vare_k, \varrho_{[k,n]}) \big) \right)
   \leq \sum_{k=1}^n m_{k,1}^2 \varrho_{[k,n]}^2
   \to 0 
 \]
 by \eqref{quad}.
Clearly
 \[
   \conv_{k=1}^n
    \CP\left( \Bi(\vare_k, \varrho_{[k,n]} ) \right)
   = \CP\left( \sum_{k=1}^n
                \Bi(\vare_k, \varrho_{[k,n]} ) \right),
 \]
 hence, in order to prove the statement, it suffices to show
 \[
   \sum_{k=1}^n \Bi(\vare_k, \varrho_{[k,n]} ) \to \mu.
 \]
We will check
 \begin{equation} \label{lim_mu}
   \sum_{k=1}^n
    \PP\left( \Bi\left(\vare_k, \varrho_{[k,n]} \right) = j \right)
   \to \mu\{j\} \qquad
   \text{for all \ $j\in\NN$.}
 \end{equation}
First note that by Taylor's formula, for all \ $p \in [0,1]$ \ and all
 \ $K,I\in\NN$,
 \[
   (1-p)^K = \sum_{i=0}^{I-1} \binom{K}{i} (-1)^i p^i + R_{K,I}(p),
 \]
 where 
 \[
   | R_{K,I}(p) | \leq \binom{K}{I} p^I.
 \]
Hence
 \begin{align*}
  \PP\left( \Bi(\vare_k, \varrho_{[k,n]} ) = j \right)
  &= \sum_{\ell=j}^\infty
      \binom{\ell}{j} \varrho_{[k,n]}^j ( 1 - \varrho_{[k,n]})^{\ell-j} \,
      \PP(\vare_k = \ell) \\
  & = \sum_{\ell=j}^\infty
       \binom{\ell}{j} \varrho_{[k,n]}^j \PP(\vare_k = \ell)
       \left[ \sum_{i=0}^{J-j-1} \binom{\ell-j}{i} (-1)^i \varrho_{[k,n]}^i
              +  R_{\ell-j,J-j}(\varrho_{[k,n]}) \right] \\
  &= \sum_{i=0}^{J-j-1} \sum_{\ell=j+i}^\infty
     \frac{(-1)^i \, \ell!}{j! \, i! \, (\ell-j-i)!} \, \varrho_{[k,n]}^{j+i} \,
     \PP(\vare_k = \ell) 
     + \tR_{n,k,j} \\
  &= \frac{1}{j!}
     \sum_{i=0}^{J-j-1} \frac{(-1)^i}{i!} m_{k,j+i} \, \varrho_{[k,n]}^{j+i}
     + \tR_{n,k,j},
 \end{align*}
 where the sum is 0 if \ $j\geq J$ \ and
 \[
   \tR_{n,k,j} := \sum_{\ell=j}^\infty
                 \binom{\ell}{j} \varrho_{[k,n]}^j \PP(\vare_k = \ell)
                 R_{\ell-j,J-j}(\varrho_{[k,n]}).
 \]
Assumption (ii) implies \eqref{lim_m} again and we have
\[
     \sum_{k=1}^n \frac{1}{j!}\sum_{i=0}^{J-j-1} \frac{(-1)^i}{i!}
     m_{k,j+i} \, \varrho_{[k,n]}^{j+i} = \frac{1}{j!}
     \sum_{i=0}^{J-j-1} \frac{(-1)^i}{i!} \sum_{k=1}^n
      m_{k,j+i} \, \varrho_{[k,n]}^{j+i} \to  \frac{1}{j!}
     \sum_{i=0}^{J-j-1} \frac{(-1)^i}{i!} \lambda_{j+i}=\mu\{j\}
\]
as \ $n\to\infty$ \ for \ $j=1,\ldots,J-1$. \ Moreover, \eqref{lim_m} 
implies
 \[
   \sum_{k=1}^n | \tR_{n,k,j} |
   \leq \sum_{k=1}^n
         \sum_{\ell=j}^\infty
          \binom{\ell}{j} \varrho_{[k,n]}^j \PP(\vare_k = \ell)
          \binom{\ell}{J-j} \varrho_{[k,n]}^{J-j}
   = \frac{1}{j! (J-j)!} \sum_{k=1}^n m_{k,J} \varrho_{[k,n]}^J
   \to \mu\{J\}=0,   
 \]
 hence we conclude \eqref{lim_mu}.
\proofend

Next we study the case when the intensity measure \ $\mu$ \ of the limiting
 compound Poisson distribution \ $\CP(\mu)$ \ may have unbounded support.

\begin{Thm}\label{CP_unbounded} \
Let \ $(X_n)_{n\in\ZZ_+}$ \ be an inhomogeneous INAR(1) process.
Assume that
 \begin{enumerate}
  \item[\textup{(i)}] \ 
   $\varr_n < 1$ \ for all \ $n \in \NN$,
    \ $\lim\limits_{n\to \infty} \varr_n = 1$,
    \ $\sum\limits_{n=1}^{\infty} (1-\varr_n) = \infty$,
  \item[\textup{(ii)}] \
   $\lim\limits_{n\to\infty} \frac{m_{n,j}}{j(1-\varr_n)}
    = \lambda_j \in [0,\infty)$
   \ for all \ $j \in \NN$ \ such that the limits
   \begin{equation}\label{mudef_infty}
    \mu\{j\} := \frac{1}{j!}
                \sum_{i=0}^\infty \frac{(-1)^i}{i!} \lambda_{j+i},
    \qquad j \in \NN,
   \end{equation}
   exist.   
 \end{enumerate}
Then
 \[
     X_n \distr \CP(\mu) 
   \quad \text{as} \quad n\to\infty.
 \]
\end{Thm}

\noindent
\textbf{Proof.}
We follow the second proof of Theorem \ref{CP_bounded} by Poisson
 approximation.
We have to show that \ $\mu$ \ is a finite measure on \ $\NN$ \ and to check
 \eqref{lim_mu}.
First note that by Taylor's formula, for all \ $p \in [0,1]$ \ and all
 \ $K,I\in\NN$,
 \[
   \sum_{i=0}^{2I-1} \binom{K}{i} (-1)^i p^i
   \leq (1-p)^K 
   \leq \sum_{i=0}^{2I} \binom{K}{i} (-1)^i p^i .
 \]
Hence for all \ $I\in\NN$,
 \begin{align*}
  \PP\left( \Bi(\vare_k, \varrho_{[k,n]} ) = j \right)
  &= \sum_{\ell=j}^\infty
      \binom{\ell}{j} \varrho_{[k,n]}^j ( 1 - \varrho_{[k,n]})^{\ell-j} \,
      \PP(\vare_k = \ell) \\
  &\leq \sum_{\ell=j}^\infty
         \binom{\ell}{j} \varrho_{[k,n]}^j \PP(\vare_k = \ell)
         \sum_{i=0}^{2I} \binom{\ell-j}{i} (-1)^i \varrho_{[k,n]}^i \\
  &= \sum_{i=0}^{2I} \sum_{\ell=j}^\infty
     \frac{(-1)^i \, \ell!}{j! \, i! \, (\ell-j-i)!} \, \varrho_{[k,n]}^{j+i} \,
     \PP(\vare_k = \ell) \\
  &= \frac{1}{j!}
     \sum_{i=0}^{2I} \frac{(-1)^i}{i!} m_{k,j+i} \, \varrho_{[k,n]}^{j+i} .
 \end{align*}
One can easily check that \eqref{lim_m} holds for all \ $j\in\NN$, \ and
we obtain
 \[
   \limsup_{n\to\infty}
    \sum_{k=1}^n
     \PP\left( \Bi(\vare_k, \varrho_{[k,n]} ) = j \right)
   \leq \frac{1}{j!}
     \sum_{i=0}^{2I} \frac{(-1)^i}{i!} \lambda_{j+i} .
 \]
In a similar way, for all \ $I\in\NN$,
 \[
   \liminf_{n\to\infty}
    \sum_{k=1}^n
     \PP\left( \Bi(\vare_k, \varrho_{[k,n]} ) = j \right)
   \geq \frac{1}{j!}
     \sum_{i=0}^{2I-1} \frac{(-1)^i}{i!} \lambda_{j+i} ,
 \]
 hence by the existence of the limits \eqref{mudef_infty} we conclude
 \eqref{lim_mu}.

Finally, for all \ $J\in\NN$, \ we have
 \begin{align*}
  \sum_{j=1}^J \mu\{j\}
  &= \sum_{j=1}^J
      \lim_{n\to\infty}
       \sum_{k=1}^n   
        \PP\left( \Bi(\vare_k, \varrho_{[k,n]} ) = j \right)
   = \lim_{n\to\infty}
      \sum_{k=1}^n  
       \PP\left( 1\leq \Bi(\vare_k, \varrho_{[k,n]} ) \leq J \right) \\
  &\leq \lim_{n\to\infty}
         \sum_{k=1}^n  
          \PP\left( \Bi(\vare_k, \varrho_{[k,n]} ) \geq 1 \right)
   \leq \lim_{n\to\infty}
         \sum_{k=1}^n 
          m_{k,1} \varrho_{[k,n]}
   = \lambda_1
 \end{align*}
 using again \eqref{Bi}.
Consequently, \ $\sum_{j=1}^\infty \mu\{j\} \leq \lambda_1 < \infty$, \ hence the
 measure \ $\mu$ \ is finite.  
\proofend

\begin{Rem}
A possible limit measure \ $\CP(\mu)$ \ in Theorem \ref{CP_unbounded} is a
 special compound Poisson measure, since its intensity measure \ $\mu$ \ has
 finite moments.
Indeed, for all \ $J,\ell\in\NN$, \ we have
 \begin{align*}
  \sum_{j=\ell}^J j(j-1)\cdots(j-\ell+1) \mu\{j\}
  &= \lim_{n\to\infty}
      \sum_{k=1}^n  
       \sum_{j=\ell}^J
        j(j-1)\cdots(j-\ell+1)
        \PP( \Bi\left(\vare_k, \varrho_{[k,n]} ) = j \right) \\
  &\leq \lim_{n\to\infty}
         \sum_{k=1}^n  
          \sum_{j=\ell}^\infty
           j(j-1)\cdots(j-\ell+1)
           \PP( \Bi\left(\vare_k, \varrho_{[k,n]} ) = j \right).
 \end{align*}
It is easy to check that for all \ $p\in [0,1]$ \ and for all non-negative
 integer-valued random variable \ $\vare$ \ we have
 \[
   \sum_{j=\ell}^\infty
    j(j-1)\cdots(j-\ell+1) \,
    \PP( \Bi\left(\vare, p \right) = j )
   = p^\ell \, \EE \vare ( \vare - 1 ) \cdots ( \vare - \ell +1 ).
 \]
Consequently,
 \[
   \sum_{j=\ell}^J j(j-1)\cdots(j-\ell+1) \mu\{j\}
   \leq \lim_{n\to\infty}
         \sum_{k=1}^n
          m_{k,\ell} \varrho_{[k,n]}^\ell
   = \lambda_\ell
   < \infty
 \]
using again \eqref{lim_m}.
\end{Rem}

\begin{Ex}
For \ $n\in\NN$, \ let \ $\varr_n=1-\frac{1}{n}$ \ and \
$\PP(\vare_n=j)=\frac{1}{nj(j+1)}$, \ $j\in\NN$, \ $\PP(\vare_n=0)
= 1-\frac{1}{n}$. \ Then \ $\EE\vare_n=\infty$ \ for all \ $n\in\NN$,
\ thus inequality \eqref{Bi_CP} is not enough to prove the compound Poisson 
convergence. Moreover, \ $\frac{\PP(\vare_n\ge j)}{j(1-\varr_n)}=\frac{1}{j^2}$
\ for all \ $j,n\in\NN$. \ The measure \ $\mu$ \ on \ $\NN$ \ defined
by \ $\mu\{j\}:=\frac{1}{j^2}$, \ $j\in\NN$, \ is finite and the infinite
series \ $\sum_{j=1}^\infty j\mu\{j\}$ \ diverges. We prove that \ $(X_n)_{n\in\ZZ_+}$ 
\ converges to \ $\CP(\mu)$ \ in spite of the fact that assumption (ii)
of Theorem \ref{CP_unbounded} does not hold. We have, for \ $n\in\NN$ \ and 
\ $z\in\CC$ \ with \ $|z|<1$ \ and \ $z\neq 0$,
\[
     H_n(z) = 1-\frac{1}{n}+\frac{1}{n}\sum_{j=1}^\infty \frac{z^j}{j(j+1)}
            = 1-\frac{1-z}{n}\left(1+\sum_{j=1}^\infty \frac{z^j}{j+1}\right)
            = 1+\frac{(1-z)\ln(1-z)}{nz},
\]
which representation is valid on the whole \ $D$. \ By Lemma \ref{Fn} we have
\[
      F_n(z) = \prod_{k=1}^n\left(1+\frac{(1-z)\ln(\frac{k}{n}(1-z))}
      {n(1-\frac{k}{n}(1-z))}\right), \qquad
      z \in D, \quad n \in \NN.
\]
Consider the functions \ $\tF_n:D\to D$, \ $n\in\NN$, \ defined by
\[
     \tF_n(z) := \prod_{k=1}^n \ee^{\frac{(1-z)\ln(\frac{k}{n}(1-z))}
      {n(1-\frac{k}{n}(1-z))}}.
\]
We have
\begin{equation}\label{tFconv}
     \tF_n(z)=\exp\left\{\frac{1}{n}\sum_{k=1}^n \frac{(1-z)\ln(\frac{k}{n}(1-z))}
      {1-\frac{k}{n}(1-z)}\right\} \to \exp\left\{(1-z)\int_0^1 
      \frac{\ln(t(1-z))}{1-t(1-z)}\ \dd t \right\}
\end{equation} 
as \ $n\to\infty$. \ Since for the dilogarithm, see  Abramowitz and Stegun
\cite[Section 27.7]{AS},
\[
      \mathrm{Li}_2(z) := \sum_{j=1}^\infty \frac{z^j}{j^2} = - \int_0^z 
      \frac{\ln(1-u)}{u}\ \dd u, \qquad z\in D,
\]
holds, we have
\[
     \tF_n(z)\to \exp\left\{\sum_{j=1}^\infty \frac{z^j-1}{j^2}\right\}\qquad
     \text{as} \quad n\to\infty
\]
for all \ $z\in D$. \ On the other hand, one can easily check that 
\begin{equation}\label{maxkn}
    \max_{1\sleq k\sleq n} \left|\frac{\ln(\frac{k}{n}(1-z))}{n(1-\frac{k}{n}
    (1-z))}\right|\to 0 \qquad \text{as \ $n \to \infty$ \ for all \ $z 
    \in D$ \ with \ $z \ne 1$.}
\end{equation}
Namely, for all \ $z \in D$ \ with \ $z \ne 1$, \ all \ $n \geq 2$ \ and all
 \ $1 \leq k\leq n$ \ we have
\[
  \left| 1 - \frac{k}{n} (1-z) \right|^2
  = 1 - \frac{2\alpha k}{n} \left( 1 - \frac{k}{n} \right) 
     - \frac{k^2}{n^2} ( 1 - |z|^2 )\leq 1 - \frac{2\alpha (n-1)}{n^2} 
   \leq  1 - \frac{\alpha}{n} < 1,     
\]
where \ $\alpha := 1 - \Re z \in (0,2]$. \ Moreover,
 \[
   \ln (1-u) = - \sum_{j=1}^\infty \frac{u^j}{j}, \qquad
   \text{for all \ $u \in \CC$ \ with \ $|u| < 1$.}
 \]
Hence, for all \ $z \in D$ \ with \ $z \ne 1$, \ all \ $n \geq 2$, \ and all
\ $1 \leq k\leq n$ \ we conclude
 \begin{align*}
  \left| \frac{ \ln \left( \frac{k}{n} (1-z) \right) }
              { n \left( 1 - \frac{k}{n} (1-z) \right)} \right|
  &\leq \frac{1}{n}
       \sum_{j=1}^\infty \frac{1}{j}
        \left| 1 - \frac{k}{n} (1-z) \right|^{j-1} \\
  &\leq \frac{1}{n}
       \sum_{j=1}^\infty \frac{1}{j}
        \left( 1 - \frac{\alpha}{n} \right)^{(j-1) / 2}
  = - \frac{ \ln \left( 1 - \sqrt{1 - \frac{\alpha}{n}} \right) }
           { n \sqrt{1 - \frac{\alpha}{n}}}
  \to 0
 \end{align*}
 as \ $n \to \infty$. \ An application of the inequality \ $|\ee^u -1-u| 
\leq |u|^2$ \ valid for all \ $u \in \CC$ \ with \ $|u| \leq 1/2$ \ implies
\[
      | \tF_n(z) - F_n(z) | \leq \sum_{k=1}^n \left|\frac{(1-z)\ln(\frac{k}{n}
      (1-z))}{n(1-\frac{k}{n}(1-z))}\right|^2 \to 0 \qquad \text{as} \quad n\to\infty
\]
for all \ $z \in D$ \ by \eqref{tFconv} and \eqref{maxkn}. Thus we finished
the proof.
\end{Ex}

\textbf{Open Problem.} The above example shows that in Theorem \ref{CP_unbounded}
we do not exhaust the possible limiting compound Poisson distribution. We conjecture 
that every compound Poisson measure can appear as a limiting distribution
of an inhomogeneous INAR(1) process. 

Theorem \ref{triangular} can also be extended for the case of limiting
compound Poisson distribution. 

\begin{Thm}
Let \ $k_n\in\NN$ \ for all \ $n\in\NN$,\ and \ $\{\zeta_{n,k}:1\leq k\leq k_n,
\  n\in\NN\}$ \ be non-negative integer-valued random variables with factorial 
moments
 \[
   m_{n,k,j} := \EE \zeta_{n,k} (\zeta_{n,k} - 1 ) \cdots (\zeta_{n,k} - j + 1 ),
   \qquad j\in\NN.
 \]
Moreover, let \ $p_{n,k}\in [0,1]$, \ $1\leq k\leq k_n$,
\ $n\in\NN$. \ Assume that
 \begin{enumerate}
  \item[\textup{(i)}] \ 
  $\sum\limits_{k=1}^{k_n} p^j_{n,k}m_{n,k,j}\to \lambda_j\in [0,\infty)$
   \ for all \ $j \in \NN$ \ such that the limits in \eqref{mudef_infty}
   exist;
  \item[\textup{(ii)}] \
   $\sum\limits_{k=1}^{k_n} (p_{n,k}\EE\zeta_{n,k})^2\to 0$;
 \end{enumerate}
as \ $n\to\infty$. \ Then
 \[
    \conv_{k=1}^{k_n} \Bi\big( \cL(\zeta_{n,k}), p_{n,k} \big)\to \CP(\mu) 
 \]
in law as \ $n\to\infty$.
\end{Thm}

\section{Appendix}

\begin{Lem}\label{compprod} \
If \ $a_k,b_k\in D$, \ $k=1,\dots,n$, \ then
 \[
   \left| \prod_{k=1}^n a_k - \prod_{k=1}^n b_k \right|
   \leq \sum_{k=1}^n | a_k - b_k |.
 \] 
\end{Lem}

\noindent
\textbf{Proof.} \
The statement follows from
 \[
   \prod_{k=1}^n a_k - \prod_{k=1}^n b_k
   = \sum_{k=1}^n
      \left( \prod_{j=1}^{k-1} a_j \right)
      (a_k - b_k)
      \left( \prod_{j=k+1}^n b_j \right)
 \]
 valid for arbitrary \ $a_k,b_k\in\CC$, \ $k=1,\dots,n$.
\proofend 

\begin{Lem}\label{conv} \
If \ $\mu_k$, $\nu_k$, \ $k=1,\dots,n$, \ are probability measures on \ $\ZZ_+$
 \ then
 \[
   d\left(\,\conv_{k=1}^n \mu_k, \, \conv_{k=1}^n \nu_k \right)
   \leq \sum_{k=1}^n d( \mu_k,\, \nu_k ).
 \] 
\end{Lem}

\noindent
\textbf{Proof.} \
The inequality
 \[
   d\left(\,\conv_{k=1}^n \mu_k, \, \conv_{k=1}^n \nu_k \right)
   \leq d\bigg( \prod_{k=1}^n \mu_k, \, \prod_{k=1}^n \nu_k \bigg)
 \] 
 easily follows from the definition of the total variation distance, where
 \ $\prod$ \ denotes product of measures.
By Barbour et al. \cite[Proposition A.1.1]{BHJ}, we have
 \[
   d\bigg( \prod_{k=1}^n \mu_k, \, \prod_{k=1}^n \nu_k \bigg)
   \leq \sum_{k=1}^n d( \mu_k,\, \nu_k ),
 \] 
 and we obtain the statement.
\proofend 

In the proofs we use extensively the following lemma about some summability
 methods defined by the sequence \ $(\varr_n)_{n\in\NN}$ \ of the offspring
 means.

\begin{Lem}\label{Toeplitz} \
Let \ $(\varr_n)_{n\in\NN}$ \ be a sequence of real numbers such that
 \ $\varr_n\in[0,1)$ \ for all \ $n\in\NN$,
 \ $\lim\limits_{n\to \infty}\varr_n=1$, \ and
 \ $\sum\limits_{n=1}^{\infty}(1-\varr_n)=\infty$.
\ Put
 \[
   a^{(k)}_{n,j} := ( 1 - \varr_j ) \prod_{\ell=j+1}^n \varr^k_\ell \qquad
   \text{for \ $n,j,k\in\NN$ \ with \ $j\leq n$.} 
 \] 
Then \ $a^{(k)}_{n,j}\leq a^{(1)}_{n,j}$ \ for all \ $n,j,k\in \NN$ 
\ with \ $j\leq n$,
 \begin{equation}\label{max_rho}
   \max_{1\sleq j\sleq n} a^{(1)}_{n,j} \to 0
   \qquad \text{as \ $n \to\infty$,}
 \end{equation}
 and for an arbitrary sequence \ $(x_n)_{n\in\NN}$ \ of real numbers with
 \ $\lim\limits_{n\to \infty}x_n=x\in\RR$, 
 \begin{equation}\label{T}
   \sum_{j=1}^n a^{(k)}_{n,j} \, x_j  \to \frac{x}{k}
   \qquad \text{as \ $n \to\infty$ \ for all \ $k\in\NN$.}
 \end{equation}
\end{Lem}

\noindent
\textbf{Proof.} \
For each \ $J\in\NN$, \ we have the inequality
 \[
   0 \leq \max_{1\sleq j\sleq n} a^{(1)}_{n,j}
     \leq \max_{j>J} a^{(1)}_{n,j} + \max_{1\sleq j\sleq J} a^{(1)}_{n,j}
     \leq \max_{j>J} ( 1 - \varr_j )
          + \max_{1\sleq j\sleq J}
             ( 1 - \varr_j ) \prod_{\ell=J+1}^n \varr_\ell,
 \]
 hence letting \ $n\to\infty$, \ we obtain
 \[
   0 \leq \limsup_{n\to\infty}\max_{1\sleq j\sleq n} 
     a^{(1)}_{n,j} \leq \max_{j>J} ( 1 - \varr_j ),
 \] 
 since
 \begin{equation}\label{prod}
   0 \leq \prod_{\ell=J+1}^n \varr_\ell
     \leq \exp \bigg\{ - \sum_{\ell=J+1}^n (1-\varr_\ell) \bigg\}
     \to 0 \qquad \text{as \ $n\to\infty$.}
 \end{equation}
Now letting \ $J\to\infty$ \ we get
 \[
   0 \leq \limsup_{n\to\infty}\max_{1\sleq j\sleq n} 
      a^{(1)}_{n,j} \leq 0,
 \] 
 and we conclude \eqref{max_rho}.

By the Toeplitz theorem, in order to prove \eqref{T}, we have to show
 \begin{align}
  \lim_{n\to\infty} a^{(k)}_{n,j} &= 0 \qquad \text{for all \ $j\in\NN$,}
   \label{T1} \\
  \lim_{n\to\infty} \sum_{j=1}^n a^{(k)}_{n,j} &= \frac{1}{k}, \label{T2} \\
  \sup_{n\in\NN} \sum_{j=1}^n |a^{(k)}_{n,j}| &< \infty \label{T3} 
 \end{align}
 for all \ $k\in\NN$.
\ By the assumptions,
 \[
   0\leq a^{(k)}_{n,j}
    \leq (1-\varr_j) \prod_{\ell=j+1}^n \varr_\ell
    \leq (1-\varr_j) \exp\bigg\{-\sum_{\ell=j+1}^n (1-\varr_\ell) \bigg\}
    \to0 \qquad \text{as \ $n\to\infty$,}
 \] 
 and we obtain \eqref{T1}.
Next we prove \eqref{T2} and \eqref{T3} for \ $k=1$.
\ We have 
 \[
   \sum_{j=1}^n a^{(1)}_{n,j}
   = \sum_{j=1}^n (1-\varr_j) \prod_{\ell=j+1}^n \varr_\ell 
   = \sum_{j=1}^n
      \left( \prod_{\ell=j+1}^n \varr_\ell
             - \prod_{\ell=j}^n \varr_\ell \right) 
   = 1 - \prod_{\ell=1}^n \varr_\ell
   \to1 \qquad \text{as \ $n\to\infty$}
 \]
 by \eqref{prod}, and
 \ $\lim\limits_{n\to\infty} \sum\limits_{j=1}^n a^{(1)}_{n,j} = 1$ \ also
 implies that
 \ $\sup\limits_{n\in\NN} \sum\limits_{j=1}^n |a^{(1)}_{n,j}|
    = \sup\limits_{n\in\NN} \sum\limits_{j=1}^n a^{(1)}_{n,j} < \infty$.
\ Hence we finished the proof of the statement of the lemma in case \ $k=1$.

The aim of the following discussion is to show \eqref{T2} and \eqref{T3} for
 all \ $k\geq2$.
\ Observe that
 \begin{align*}
  1 - \prod_{\ell=1}^n \varr^k_\ell
  &= \sum_{j=1}^n
      \left( \prod_{\ell=j+1}^n \varr^k_\ell
             - \prod_{\ell=j}^n \varr^k_\ell \right)
   = \sum_{j=1}^n (1-\varr^k_j) \prod_{\ell=j+1}^n \varr^k_\ell \\
  &= \sum_{j=1}^n
      \bigg( \sum_{i=1}^k \binom{k}{i} (1-\varr_j)^i \bigg)
      \prod_{\ell=j+1}^n \varr^k_\ell \\
  &= k \sum_{j=1}^n a^{(k)}_{n,j}
     + \sum_{i=2}^k
        \binom{k}{i} \sum_{j=1}^n (1-\varr_j)^{i-1} a^{(k)}_{n,j},
 \end{align*}
 where, by \eqref{prod},
 \ $0\leq\lim\limits_{n\to\infty} \prod\limits_{\ell=1}^n \varr^k_\ell
     \leq\lim\limits_{n\to\infty} \prod\limits_{\ell=1}^n \varr_\ell = 0$.
\ Moreover,
 \[
   0 \leq \sum_{j=1}^n (1-\varr_j)^{i-1} a^{(k)}_{n,j}
      \leq \sum_{j=1}^n (1-\varr_j)^{i-1} a^{(1)}_{n,j} \to 0
  \qquad \text{as \ $n\to\infty$ \ for all \ $i\geq2$} 
 \]
 by the lemma for \ $k=1$ \ and by the assumption
 \ $\lim\limits_{n\to \infty}\varr_n=1$.
\ Consequently, we obtain \eqref{T2} and hence \eqref{T3} for all \ $k\geq2$.
\proofend 

\begin{Lem}\label{Taylor} \
Let \ $\vare$ \ be a nonnegative integer-valued random variable with factorial
 moments
 \[
   m_k := \EE \vare ( \vare - 1 ) \cdots ( \vare - k + 1 ),
   \qquad k\in\NN,
 \]
$m_0:=1$, \  and with generating function \ $H(z) = \EE(z^\vare)$, 
\ defined for \ $z \in D$. 
\ If \ $m_k<\infty$ \ for some \ $k\in\NN$ \ then 
 \[
   H(z) = \sum_{j=0}^{k-1} \frac{m_j}{j!} (z-1)^j + R_k(z)
   \qquad \text{for all \ $z \in D$,}
 \]
 where
 \[
   |R_k(z)| \leq \frac{m_k}{k!} |z-1|^k
   \qquad \text{for all \ $z \in D$.}
 \]
\end{Lem}

\noindent
\textbf{Proof.} \
By \ $m_j = \sum\limits_{\ell=0}^\infty
             \ell(\ell-1)\cdot(\ell-j+1) \, \PP(\vare=\ell)$,
 \[
   R_k(z)
   = H (z) - \sum_{j=0}^{k-1} \frac{m_j}{j!} (z-1)^j
   = \sum_{\ell=0}^\infty
      \bigg( z^\ell -  \sum_{j=0}^{k-1} \binom{\ell}{j} (z-1)^j \bigg)
      \PP(\vare =\ell),
 \]
 and by Taylor's formula for the function \ $z\mapsto z^\ell$ \ we get
 \begin{align*}
   \bigg| z^\ell - \sum_{j=0}^{k-1} \binom{\ell}{j} (z-1)^j \bigg|
   &\leq \frac{1}{k!}
         |z-1|^k
         \sup_{\theta\in[0,1]}
          \big|\ell(\ell-1)\cdots(\ell-k+1)(1+\theta(z-1))^{\ell-k}\big| \\
   &\leq \frac{1}{k!}
         \ell(\ell-1)\cdots(\ell-k+1)
         |z-1|^k
 \end{align*}
for all \ $z \in D$.
\proofend

\vspace*{10mm}    

\noindent
\textsc{L\'aszl\'o Gy\"orfi} and \textsc{Katalin Varga}, \\
Department of Computer Science and Information Theory, \\
Budapest University of Technology and Economics, \\
Stoczek u. 2, Budapest, Hungary, H-1521; \\
\textit{e-mails}: {\tt \{gyorfi,varga\}@szit.bme.hu}\\[5mm]
\textsc{M\'arton Isp\'any} and \textsc{Gyula Pap}, \\
Department of Applied Mathematics and Probability Theory, \\
Faculty of Informatics, University of Debrecen, \\
Pf.12, Debrecen, Hungary, H-4010;\\
\textit{e-mails}: {\tt \{ispany,papgy\}@inf.unideb.hu}

\end{document}